\documentclass[english]{amsart}
\usepackage[T1]{fontenc}
\usepackage{amsthm}
\usepackage{amstext}
\usepackage{amssymb}

\makeatletter
\numberwithin{equation}{section}
\numberwithin{figure}{section}
\theoremstyle{plain}
\newtheorem{thm}{\protect\theoremname}

\makeatother

\usepackage{babel}
\providecommand{\theoremname}{Theorem}

\begin{document}

\title{Corrections to \textquotedbl{}On continued fractions of given period\textquotedbl{}}

\author{Vladimir Pletser}

\address{European Space Research and Technology Centre \\
ESA-ESTEC P.O. Box 299 \\
NL-2200 AG Noordwijk \\
The Netherlands}

\address{European Space Research and Technology Centre \\
ESA-ESTEC P.O. Box 299 \\
NL-2200 AG Noordwijk \\
The Netherlands}

\email{Vladimir.Pletser@esa.int}


\begin{abstract}
Corrections are brought to an article of Friesen on continued fractions
of a given period. 
\end{abstract}
\maketitle
Friesen has proven \cite{1} that, for any $k\in\mathbb{Z^{\dotplus}}$,
there are infinitely many squarefree integers $N$, where the continued
fraction expansion of $\sqrt{N}$ has period equal to $k$. This was
demonstrated in a Corollary following a Theorem stating that
\begin{thm}
Let $\left\lfloor \sqrt{N}\right\rfloor $ denote the greatest integer
$\leq\sqrt{N}$. Then the equation  
\[
\sqrt{N}=\left[\left\lfloor \sqrt{N}\right\rfloor ;\overline{a_{1},a_{2},\text{\ldots},a_{k-2}=a_{2},a_{k-1}=a_{1},a_{k}=2\left\lfloor \sqrt{N}\right\rfloor }\right]
\]

has, for any symmetric set of positive integers $\left\{ a_{1},...,a_{k-1}\right\} $,
infinitely many squarefree solutions $N$ whenever either $Q_{k-2}$
or $\left(Q_{k-2}^{2}-\left(-1\right)^{k}\right)/Q_{k-1}$ is even.
If both quantities are odd, then there are no solutions $N$ even
if the squarefree condition is dropped.
\end{thm}
where $Q_{k}$ are convergents of the of the continued fraction of
$\sqrt{N}$. This theorem is demonstrated using three Lemmas. 

Unfortunately, this paper contains two mistakes, that we want to correct
here.

First, in the demonstration of the Lemma 1, (see \cite{1}, p. $12$,
line 19), the expression of $N$ should be $N=N\left(b\right)=\alpha b^{2}+\beta b+\gamma$
(instead of $\alpha b^{2}+\beta b^{2}+\gamma$, which would not make
sense). However this does not change the final result of this Lemma
1.

Second, Friesen proved the Corollary with a sufficient condition that
for each $k\in\mathbb{Z^{\dotplus}}$, a symmetric set of positive
integers $\left\{ a_{1},...,a_{k-1}\right\} $ exists such that $Q_{k-1}$
is odd. The case $k=1$ is direct as it yields $N(b)=b^{2}+1$ giving
an infinity of squarefree $N$. 

However, for the case $k>1$, the demonstration (see \cite{1}, p.
$13$, lines 19-25) is wrong as it contradicts the final statement
made in the conclusion (see \cite{1}, p. $13$, lines 35-36).

I have attempted to correct his demonstration herebelow. 

\textquotedbl{}Assume $k>1$. If $k\equiv0\left(mod\,3\right)$ {[}instead
of $k\neq0\left(mod\,3\right)${]} then set $a_{1}=a_{k-1}=2$ and
$a_{i}=1$ for $i=2,...,k-2$. If $k\neq0\left(mod\,3\right)$ {[}instead
of $k\equiv0\left(mod\,3\right)${]} set$a_{i}=1$ for $i=1,...,k-1$.
In both cases we have the recursion formula for $Q_{n}$ giving us
copies of the Fibonacci sequence (...). In the first instance {[}i.e.
$k\equiv0\left(mod\,3\right)${]} we have $Q_{n}=F_{n+2}$for $n=2,...,k-2$
{[}instead of $i=1,...,k-2${]} and $Q_{k-1}=2Q_{k-2}+Q_{k-3}=2F_{k}+F_{k-1}$
{[}instead of $2F_{k}+F_{k-1}+F_{k+2}${]}. If $k\neq0\left(mod\,3\right)$
we have $Q_{n}=F_{n+1}$ for $n=1,...,k-1$, hence $Q_{k-1}=F_{k}$.
But, as $F_{k}$ is even only when $k\equiv0\left(mod\,3\right)$,
we see that $Q_{k-1}$ is odd in either situation. Therefore, either
$Q_{k-2}$ or $\left(Q_{k-2}^{2}-\left(-1\right)^{k}\right)/Q_{k-1}$
is even and we have satisfied the conditions of the Theorem, thus
proving the Corollary.\textquotedbl{}

With these corrections, these statements agree with those of the conclusion
namely ``By setting $a_{i}=1$ for $i=1,...,k-1$ if $k\neq0\left(mod\,3\right)$
(and $a_{1}=a_{k-1}=2$, $a_{i}=1$ for $i=2,...,k-2$ if $k\equiv0\left(mod\,3\right)$)
it was shown that the conditions of the Theorem are met.'' 

Despite these mistakes, the final result still stands.

\end{document}